\documentclass[10pt]{article}

\usepackage[all]{xy}

\usepackage{hyperref}

\usepackage{amsmath,amssymb}
\usepackage[latin1]{inputenc}
\usepackage{lscape}
\usepackage{mathtools}

\newtheorem{thm}{Theorem}[section]
\newtheorem{conj}[thm]{Conjecture}

\newtheorem{cor}[thm]{Corollary}
\newtheorem{prop}[thm]{Proposition}
\newtheorem{defi}[thm]{Definition}

\newtheorem{ques}[thm]{Question}

\def\Hom{\mathrm {Hom}}

\def\Im{{\rm Im}}

\def\ra{\longrightarrow}

\numberwithin{equation}{section}

\def\GL{\mathbb {GL}}
\def\M{\mathbb {M}}
\def\Q{\mathbb Q}
\def\Z{\mathbb Z}
\def\F{\mathbb F}

\def\F2{{\mathbb F}_2}
\def\U{\mathcal U}
\def\calS{\mathcal S}

\def\F2{{\mathbb F}_2}

\begin{document}

\title{Questions and conjectures about the modular representation theory of the general linear group $\GL_n(\F2)$ and the Poincar\'e series of unstable modules}
\author{Kirian Delamotte, Nguyen D. H. Hai, Lionel Schwartz} 
\date {UMR 7539 CNRS, 
 Universit\'e Paris 13\\
 University of Hu\'e \\
 LIA CNRS Formath Vietnam
 }

\maketitle

\begin{abstract} This note is devoted to some questions about the representation theory over the finite field $\F2$ of the general
linear groups $\GL_n(\F2)$ and Poincar\'e series of unstable modules. The first draft was describing two conjectures. They were presented during talks made at VIASM in summer 2013. Since then one conjecture has been disproved, the other one has been proved. These
results naturally  lead to new questions which are going to be discussed. 
In winter 2013,  Nguyen Dang Ho Hai proved the second conjecture, he  disproved the first one in spring 2014. Up to now, the proof of  the second one depends on a major topological result: the Segal conjecture. This discussion could be extended to an odd prime, but we will not do it here, just a small number of remarks will be made.

\end{abstract}

\section{Introduction}
Over the finite field $\F2$, the general linear group group $\GL_n(\F2)$ has, up to isomorphism, $2^{n-1}$ distinct simple representations. They are  indexed by strictly decreasing partitions $\lambda$ starting 
with the integer $n$ \cite{JK81}. These partitions  are said to be  $2$-regular. As the indexing set for simple representations of the the general linear group depends on the reference, we will recall the different equivalent choices of indexation in  ``topological references".

Let $p$ be  a prime number. A partition $\lambda$ is said to be $p$-regular if and only if it does not contain $p$ non zero successive parts which are equal. At the prime $2$ these are the strictly decreasing partitions.  Let $\lambda$ be  a partition, the associated partition $\lambda'$ is  defined by $\lambda'_j =\#\{i, \vert \lambda_i \geq j\}$.  The preceding condition on $\lambda$ translates for $\lambda'$ into: for each $i$, we have $\lambda'_{i}-\lambda'_{i+1} \leq p-1$ (a partition sharing this property  is said to be column $p$-regular). Finally, in \cite{Harris-Shank-92} the representations are indexed
by the differences $(\lambda'_1-\lambda'_2,\ldots,\lambda'_{h-1}-\lambda'_h,\lambda'_h)$.

For example, the trivial representation of $\GL_n(\F2)$ is, in our setting, indexed by $(n)$,  the associated partition is $(\underbrace{1,\ldots,1}_\textrm{$n$ times})$. The trivial representation is indexed
by $(\underbrace{0,\ldots,0}_\textrm{$n-1$ times},1)$ in Harris and Shank's paper. The Steinberg representation is indexed by $(n,n-1, \ldots,1)$ which is equal to its associated partition while in Harris and Shank's paper the representation is indexed by 
$(\underbrace{1,\ldots,1}_\textrm{$n$ times})$ .

There is in \cite{JK81} a theoretical description, of combinatorial nature, for the simple representations. This is done using row and column stabilizers of the associated Young diagrams.  For a partition
$\lambda$ of the integer $n$ denote by $C_\lambda$, resp. by $R_\lambda$, the column, resp. row, stabilizer of the associated Young diagram. Denote by $\bar C_\lambda$, resp. $\bar R_\lambda$, the sum of their elements in the group algebra over $\F2$ of the symmetric group $\calS_n$. Let $V=(\Z/2)^{n}$, the simple representation $S_\lambda$ is isomorphic to $V^{\otimes n}\bar C_\lambda\bar R_\lambda\bar C_\lambda$ (see \cite{JK81} chapter 7 and 8, and \cite{PS98}). At an odd prime the definition of $\bar C_\lambda$ is $\sum_\sigma \mathrm{sgn}(\sigma)[\sigma]$.

\footnote{partially supported by the program ARCUS Vietnam of the french MAE and the R\'egion Ile de France}

However, the modular 
character of these representations, as well as their projective covers and the Cartan matrices, are  unknown. If $\varrho$ is a simple representation, indexed by $\lambda$, it will be  denoted $S_\lambda$;
its projective cover will be  denoted by $P_\varrho$ or  $P_\lambda$.

Let us give more examples. Let $n=2$, $V = (\Z/2)^2$,  the simple representations are indexed by $(2)$, $(2,1)$. They  are the trivial representation $\F2$ and the standard one $V$. The second one is projective; the first one has, as projective cover, a non-trivial extension by itself. Next, let $n=3$ and
$V \cong (\Z/2)^3$, the simple representations are indexed by $(3)$, $(3,1)$, $(3,2)$ and $(3,2,1)$.  They are  the trivial representation $\F2$, the standard one $V$, the contragredient of the standard one $\Lambda^2(V) = (V^{-1})^*$ and the Steinberg representation. The Steinberg representation is the kernel of $\Lambda^2(V) \otimes V \ra \Lambda^3(V) \cong \F2$; it is projective.

\vskip 2mm

The $p$-regular partitions serve also as indexing set for the simple objects of the category $\mathcal F$  of functors from the category of finite dimensional $\F2$-vector spaces to the category
of all $\F2$-vector spaces. Note that in  \cite{PS98} the simple objects of $\mathcal F$  are indexed using the associated partitions. It is often better to use the ``functorial" framework which is done implicitly below.

Let $\lambda=(\lambda_1,\ldots,\lambda_h)$, $\lambda_h>0$, be a strictly decreasing partition starting with $\lambda_1=n$;  $h$ is the length of the partition and is denoted by $\ell(\lambda)$. 
Here is a general proposition that can be found in the last chapter of \cite{JK81}:

\begin{prop}\label{basis} Let $V$ be the vector space $(\Z/2)^n$.
The simple representation $S_\lambda$ associated with $\lambda$ occurs one time in $\Lambda^\lambda(V) :=\Lambda^{\lambda_1}(V)\otimes \ldots \otimes \Lambda^{\lambda_h}(V).$ All other simple representations  occurring 
in the Jordan-H\"{o}lder decomposition are indexed by partitions lower than $\lambda$ in the dominance order.
\end{prop}

The  dominance order $\precsim$ on partitions is the following one. Let  $\lambda$ and $\mu$ be partitions. If $\sum_i \lambda_i< \sum_i \mu_i$ then $\lambda \precsim \mu$. Next if $\lambda$ and $\mu$ are partitions of the same integer ($\sum_i \lambda_i=\sum_i \mu_i$),
then  $\lambda \precsim \mu$ if  $\sum_{i<h} \lambda_i= \sum_{i<h} \mu_i$ for all $h$. With this definition the partition $(n)$ is the smallest one (see also \cite{JK81} page 23) among partitions of $n$. 

\begin{cor} 
The representations  $\Lambda^{\lambda_1}(V)\otimes \ldots \otimes \Lambda^{\lambda_h}(V)$,  $\lambda$ strictly decreasing partition starting with $n$,
form a basis of the Grothendieck representation ring $R_{\F2}(\GL_n(\F2))$. 
\end{cor}

What follows provides a  link between modular representation theory of the general linear groups and topology. 
The group  $\GL_n(\F2)$ acts naturally on  $B(\Z/2)^n$ and thus on the mod $2$ singular cohomology $H^*B(\Z/2)^n$. S. Mitchell and S. Priddy \cite{MP83} used this action to prove the following theorem:

\begin{thm} The space $B(\Z/2)^n$ splits up, as a  spectrum, as the following wedge:
$$
\bigvee_{\varrho \in Irr_{\F2}(\GL_n(\F2))}\, \, M(\varrho)^{\vee \dim(\varrho)}.
$$
\end{thm}

The proof depends on the following observation: an element in the group algebra of $\GL_n(\F2)$ induces a map on the suspension spectrum
of $B(\Z/2)^n$. As an example of element in the group algebra we can take an idempotent associated to the projective cover of a simple representation $\varrho$. 
The spectrum  $M(\varrho)$ is the telescope of the associated map.

\vskip 2MM

A major issue is to understand the Poincar\'e series of the cohomology $H^*M(\varrho)$.
\vskip 2mm
A second link between representation theory and topology  is provided by Lannes'v $T_V$-functor. The functor $T_V$ is left adjoint to $M \mapsto M \otimes H^*V$ from the category $\U$
of unstable modules over the Steenrod algebra  to itself. It is exact and commutes with tensor products \cite{Lan92}, also $T_V \cong T^n$. Moreover \cite{LS89} shows that 
each indecomposable reduced injective unstable module is a direct summand in some $H^*V$. Here ``reduced" means that the module does not contain a non-trivial suspension \cite{Sch94}.
It is known (Adams-Gunawardena-Miller, Lannes-Zarati) that $T(H^*V) \cong \F2[V^*] \otimes H^*V$. 

\begin{defi} Denote by $K^{red}(\U)$ the Grothendieck group of reduced injective unstable modules which are finite direct sums of indecomposable ones. \end{defi}

 It follows from \cite{LZ86} that this is a ring and from
\cite{LS89} that a basis is indexed by simple representations of all $\GL_n(\F2)$, thus by strictly decreasing partitions. 

Finally the above result of  \cite{Lan92} shows that $T$ acts as a ring homomorphism. It is of interest to understand this homomorphism.

\section{Linear (in)dependance of the Poincar\'e series of the indecomposable summands of $B(\Z/2)^n$}

Recall from the introduction that there is a stable decomposition:
$$B(\Z/2)^n\simeq
\bigvee_{\varrho \in Irr_{\F2}(\GL_n(\F2)}\, \, M(\varrho)^{\vee \dim(\varrho)}.
$$
We will denote the spectrum $M(\varrho)$  by $M(\lambda)$, using the indexing partition $\lambda$ of $\varrho$. Given a partition $\lambda=(\lambda_1,\ldots,\lambda_h)$ denote by $\bar \lambda$
the partition $(\lambda_2,\ldots,\lambda_h)$. The above decomposition is refined by (see also \cite{HK88}):

\begin{thm} The spectrum $M(\lambda)$ splits up as the wedge of two indecomposable spectra. One spectrum, denoted by $L(\lambda)$ occurs as a summand of
$B(\Z/2)^n$, but not of $B(\Z/2)^h$, $h<n$. We have:
$$
M(\lambda) \simeq L(\lambda) \bigvee L(\bar \lambda).
$$
\end{thm}

We will denote by $L_\lambda$ (resp. $M_\lambda$)  the cohomology  $H^*L(\lambda)$ (resp. $H^*M(\lambda)$), the modules  $L_\lambda$ are a list of representatives of reduced injective unstable modules \cite{LS89}.

The cohomology $H^*BV$ is isomorphic to $S^*(V^*)$,  the symmetric algebra over $V^*$. It follows from the telescope construction that the Poincar\'e series, $P_\lambda(q)$, of $M(\lambda)$ is given by the following  formula:
$$
P_\lambda(q) = \sum_{k\ge 0} a_{\lambda,k}\cdot q^k,
$$
where $a_{\lambda,k}=\dim_{\F2}\Hom_{\GL_n(\F2)}(P_\varrho, S^k(V))$ is the multiplicity of the representation $\varrho=S_\lambda$ in $S^k(V)$.

 Mitchell and Priddy showed in \cite{MP83} the following:
\begin{prop}\label{pole}
We have:
$$
P_\lambda(q) = \frac{\Pi_\lambda(q)}{(1-q^{2^n-2^{n-1}})(1-q^{2^n-2^{n-2}})\ldots(1-q^{2^n-1})}
$$
for a certain polynomial $\Pi_\lambda$ with non-negative coefficients.
This series has a pole of order $n$ at $1$.
\end{prop}

This is  because $H^*B(\Z/2)^n$ is a free, finitely generated,  module over the Dickson algebra $(H^*B(\Z/2)^n)^{\GL_n}$, which is a polynomial algebra on generators in degrees $2^n-1, 2^n-2, \ldots, 2^n-2^{n-1}$. 

Carlisle and Walker \cite{CW89} established another form with better denominator for the Poincar\'e series. We have
$$
P_\lambda(q) = \frac{\Gamma_\lambda(q)}{(1-q)(1-q^{3})\ldots(1-q^{2^n-1})}
$$
for a certain polynomial $\Gamma_\lambda$.

The polynomials $\Pi_\lambda$ and $\Gamma_\lambda$ are unknown in general. They can be computed in special cases, in particular the series of the Steinberg summand $L_n:=L_{(n,n-1,\ldots,1)}$ of Mitchell and Priddy is given by:

$$
\frac{q^{1+(2^2-1)+ \ldots +(2^n-1})}{(1-q)(1-q^3)\ldots(1-q^{2^n-1})}
$$

The connectivity of $L_\lambda$ is known, \cite{FS90}:

\begin{prop}\label{connex}
The value of the first non zero coefficient of 
$\Pi_\lambda$ is $1$ and in degree $\lambda_1+ 2 \lambda_2+ \ldots + 2^{h-1}\lambda_h$.
\end{prop}

As said above $\Pi_\lambda$ and $\Gamma_\lambda$ are not known in general. Outside of the Steinberg representation the other known examples are, in some sense, close to the Steinberg representation \cite{CW89}. For an odd prime $p$  the situation is worse, the analogous result of Proposition \ref{connex} is not even known outside of particular cases \cite{MW02}. 
\vskip 2mm

It is natural to ask whether or not the Poincar\'e series of $L_\lambda$ determines $\lambda$.  The first observation is that the order of the pole at $1$ of the series is  $\lambda_1$. Next up to $n=4$ the first non zero term of the series allows to distinguish representations. This is no longer true if $n>4$, because the Poincar\'e series of 
the representations associated to $(5,4)$ and $(5,2,1)$ have the same connectivity. But the preceding results suggest  the following:

%

\begin{conj} \label{independent}The Poincar\'e series $P_\lambda$ are linearly independent.
\end{conj}

From a topological point of view this looks reasonable. However, based on work of D. Carlisle \cite{Carlisle85} and of S. Mitchell \cite{Mit85}, the second author  has given a counterexample. Indeed, the following virtual module (with some classical notation):
$$X=[M_{(3)}]-[M_{(3,2,1)}]-[M_{(4)}]-[M_{(4,1)}]-[M_{(4,3)}]+
[M_{(4,2,1)}]+[M_{(4,3,2)}]+[M_{(4,3,1)}]+2\,[M_{(4,3,2,1)}]$$
has trivial Poincar\'e series. There is also an important observation here: the Poincar\'e series of $T(X)$ is non trivial (this is proved by direct computation). This is, from a topological point of view, a curious phenomenon and suggests that an intricate phenomenon occurs. The only hint  is  the
 existence of the following non-split short exact sequence \cite{JP94}:
$$
\{0\} \ra K \ra L_{(1)}^{\otimes 4} \ra \ L_{(2,1)}\ra \{0\}.
$$



\section{Eigenvalues and eigenspaces of Lannes' $T$ functor}

Let us consider now the Grothendieck ring, $K^{red}(\U)$, of  injective,   reduced unstable modules which are finite direct sums of indecomposable ones. This ring is filtered by subgroups 
$K^{red}_n(\U)$ generated by the classes of the $L_\lambda$ such that  $\lambda_1 \leq n$. Recall that Lannes' $T$-functor, \cite{Lan92},  is left adjoint to $M \mapsto H^*\Z/2 \otimes M$; it is exact and preserves tensor products. Thus it
acts as a ring homomorphism on  $K^{red}(\U)$. Because of the Adams-Gunawardena-Miller theorem it  preserves the filtration. 
Suppose given a basis of $K^{red}(\U)$ which is the union of basis of the $K^{red}_n$. Harris and Shank \cite{Harris-Shank-92} have proved that: 

\begin{prop} The $(2^n,2^n)$ matrix $t_n$ of the restriction of $T$ to $K^{red}_n$ writes as
$$
\begin{pmatrix} t_{n-1} & \delta_{n-1}\cr  0 & \tau_n\cr \end{pmatrix}
$$
with  $\tau_n=t_{n-1}+\delta_{n-1}$.
\end{prop}

The $2^{n-1}$ first columns of the above matrix correspond to the action of $T$ on $K^{red}_{n-1}$.
\begin{conj}\label{eigen}\ 
\begin{enumerate}
\item The eigenvalues of $T$ on $K^{red}_n$ are all powers of $2$, more precisely they are  $(1)^{2^{n-1}},(2)^{2^{n-2}},\ldots, (2^{n-1})^1,  (2^{n})^1$ (multiplicity is shown as exponent);
\item $T$ is diagonalisable on $K^{red}_n$.
\end{enumerate}

\end{conj}

The analogous conjecture can be made at any prime $p$.
The conjecture has been checked up to $n=9$ by the first author  for $p=2$, and finally proved by the second one using the Segal conjecture \cite{NDHH13}.

\section{The Poincar\'e series and the image of the homomorphism $M \mapsto P_M(q)$}

This section discusses questions  about the image of the map which sends an unstable module to its Poincar\'e series: to which extent it determines the original object.  First we assume this can be defined and that any module we consider are of finite dimension in any degree. There is a first
necessary restriction, to restrict to unstable modules which have as injective hull a finite direct sum of indecomposable injective unstable modules. Otherwise we could consider infinite direct sums
$\bigoplus_N\Sigma^n \F2^{\oplus \alpha_n}$ to represent any formal power series.

\vskip 2mm

The first thing to discuss is what should be the source of this map, it is reasonable to conjecture the following:

\begin{conj} The full subcategory of $\U$ whose objects are unstable modules whose injective hull is a finite direct sum of indecomposable injective unstable modules is thick. In particular quotients of such
objects have the same property.
\end{conj}

This is supported by Steven Sam's proof of the artinian conjecture and by \cite{CS14}.
If this is true we could form the associated Grothendieck ring which we denote by $G^{fih}(\U)$. The map $M \mapsto P_M(q)$ induces a ring homomorphism
$$
\pi\colon G^{fih}(\U) \ra \Z[[q]].
$$
It has been said above that this homomorphism is not injective by restriction to the subring $K^{red}(\U)$. 
However it is injective by restriction
to the Grothendieck subring, $G(\U)$, of finitely generated unstable modules \cite{Sch06}.

Two other Grothendieck rings are useful. The first one is the ring $K(\U)$ which has the same definition as $G^{fih}(\U)$, but the finiteness condition on the injective hull is removed. Next it is natural to introduce a topology on this ring and to consider its completion $\widehat K(\U)$.

There are no obvious conjectures about the kernel of $\pi$, on the other hand the image looks to be easier to understand. Below is a list of questions.

\begin{ques} The first question is to describe $\Im(\pi) \otimes \Q \subset \Q[[q]]$. It looks reasonable to conjecture it is a polynomial algebra:
$\Q[q] \otimes \Q[\chi_{2k+1, \, \, k \geq 0}] \otimes \Q[P_\lambda , \, \, \lambda \, \, regular]/(relations) $, where $\chi_{2k+1}=\sum_{h\ge 0} q^{(2k+1)2^h}$ and $P_\lambda$ is as above.
\end{ques}

This question contains in fact two. The first one, may be the most interesting and new, is that whether or not the image is generated by $q$ the $\chi_{2k+1}$s and the $P_\lambda$s. The second one 
is to understand the last factor, that is to explicit the relations, information about this can be found in Carlisle-Kuhn \cite{CK96} and 
Harris-Shank \cite{Harris-Shank-92}.

The last question will be formulated as a conjecture. It concerns the explicit form of the Poincar\'e series $P_\lambda$, more precisely it concerns the form of the Poincar\'e series of the virtual modules
in $K^{red}(\U)$ who are eigenvectors for $T$ associated to the eigenvalue $1$.  Numerical evidences support:

\begin{conj}The Poincar\'e series of each eigenvector for $T$ associated to the eigenvalue $1$ has no pole at $1$.
\end{conj}

If true, using work of Carlisle and Walker \cite{CW89}, the Poincar\'e series of an eigenvector associated to $1$ in $K_n^{red}$ could be chosen in such a way that the Poincar\'e series is of the form:
$$
\frac{P(q)}{\prod_{1 \leq i \leq n}(1+q+q^2+q^3+\ldots+q^{2^{i}-2})},
$$
where $P(q)$ is a polynomial.

\section{How to write algorithms and numerical results}

The first point is that we will work in $R_{\F2}(\GL_n(\F2))$ with the basis given by $\Lambda^\lambda(\F2^n)$, $\lambda$  strictly decreasing partitions starting with $n$ (see Proposition \ref{basis}). We will also work with
$R_{\F2}(\M_n(\F2))$, in this case a basis is given by the $\Lambda^\lambda(\F2^n)$, $\lambda$ $2$-regular and $\lambda_1 \leq n$.

Conjectures \ref{independent} and \ref{eigen} can be tested using computer. For the first one as the answer turns out to be negative we will not make long comments. Let us just say that it is enough to understand the span of  the symmetric powers $S^k(\F2^n)$ in the Grothendieck ring $R_{\F2}(\M_n(\F2))$. To do this write the decomposition of the symmetric powers on the basis, and use Gauss' method to compute the rank. However as the system is infinite it only yields an answer if the answer is positive. In order to write the decomposition  there are two possible ways. The first one is to use Koszul complex. It expresses  $S^k$ (in the representation ring) as an alternating sum of tensor products $S^h \otimes \Lambda ^{k-h}$, $h<k$, and do an iterative process. 
The other possibility is to use the Taylor series of $S^k$ \cite{Tro05}. It gives a decomposition of the symmetric powers in terms of tensor products of exterior powers
in the Grothendieck group.

Concerning the Conjecture \ref{eigen}, in order to write down an algorithm,  we need the following (\cite{Harris-Shank-92} Theorem 3.8):

\begin{thm}
The transpose of the matrix of $T$ on $K^{red}_n(\U) \subset K^{red}(\U)$ is equivalent to the matrix of   the following map $R_{\F2}(\M_n(\F2)) \ra R_{\F2}(\M_n(\F2))$:
$$
[\varrho] \mapsto  \sum_{0 \leq i \leq n}\, \, [\varrho \otimes \Lambda^i].
$$
\end{thm}

We have a similar result for the matrix of $T$ on the quotient $K^{red}_n/K^{red}_{n-1}$ by replacing the Grothendieck ring $R_{\F2}(\M_n(\F2))$ by 
$R_{\F2}(\GL_n(\F2))$. 
 
In both cases we need another ingredient. In the obvious iterative process we have to decompose tensor products 
 $(\Lambda^{\lambda_1}\otimes \ldots \otimes \Lambda^{\lambda_h}) \otimes \Lambda^k$ in the  basis.  Either $k$ is different from all $\lambda_i$ and there is nothing to do, but to reorder. If 
 $k=\lambda_i$ for some $i$ the tensor product is not a basis element. We have to write its decomposition in the basis. For this we use the following formula:
 $$
 [\Lambda^k \otimes \Lambda^k]=[\Lambda^k]+2[\Lambda^{k+1} \otimes \Lambda^{k-1}]-2[\Lambda^{k+2} \otimes \Lambda^{k-2}]+\ldots +(-1)^{i-1}2[\Lambda^{k+i} \otimes \Lambda^{k-i}]+\ldots
 $$
 This comes from the complex
$$
\{0\}\ \ra \Lambda^{2k}\ra \Lambda^{2k-1}\otimes \Lambda^1 \ra \ldots \ra \Lambda^k \otimes \Lambda^k \ra \ldots \ra \Lambda^1 \otimes \Lambda^{2k-1} \ra \Lambda^{2k} \ra \{0\}
$$
which is exact except at the middle, where the homology is $\Lambda^k$.

We iterate this process, a moment of reflexion shows that it has to stop sometimes.

\newpage

\subsection{Some relations in $R_{\F2}(\M_4(\F2))$ }
 
\

Without mod $2$ reduction, we have the following relations in $R_{\F2}(\M_4(\F2))$, with obvious notations:
{\scriptsize
\begin{eqnarray*}
s_{0} & = & 1(0)\\
s_{1} & = & (1)\\
s_{2} & = & (1)+(2)\\
s_{3} & = & (1)+2(2)+(3)\\
s_{4} & = & (1)+(2)+(2,1)+(4)\\
s_{5} & = & (1)+2(2)+(2,1)+(3,1)\\
s_{6} & = & (1)+3(2)+(2,1)+(3)+2(3,1)-2(4)+(4,1)\\
s_{7} & = & (1)+2(2)+2(2,1)+2(3,1)+(2,3)\\
s_{8} & = & (1)+3(2)+2(2,1)+4(3,1)+(2,3)-3(4)+(4,1)+(4,2)\\
s_{9} & = & (1)+4(2)+2(2,1)+(3)+5(3,1)+2(3,2)-4(4)+(4,1)+2(4,2)\\
s_{10} & = & (1)+3(2)+3(2,1)+5(3,1)+4(3,2)-2(4)-1(4,1)+3(4,2)+(4,3)\\
s_{11} & = & (1)+4(2)+3(2,1)+6(3,1)+3(3,2)+(3,2,1)-4(4)+(4,3)\\
s_{12} & = & (1)+5(2)+3(2,1)+(3)+8(3,1)+4(3,2)+(3,2,1)-5(4)+(4,2)+(4,2,1)\\
s_{13} & = & (1)+4(2)+4(2,1)+9(3,1)+6(3,2)+(3,2,1) 
           -4(4)-1(4,1)+2(4,2)+2(4,2,1)\\
s_{14} & = & (1)+5(2)+4(2,1)+10(3,1)+6(3,2)+2(3,2,1) 
           -6(4)-1(4,1)+2(4,2,1)+(4,3,1)\\
s_{15} & = & (1)+6(2)+4(2,1)+(3)+13(3,1)+7(3,2)+2(3,2,1) 
           -8(4)+2(4,2)+4(4,2,1)-2(4,3)+(4,3,1)\\
s_{16} & = & (1)+5(2)+5(2,1)+14(3,1)+10(3,2)+2(3,2,1) 
           -9(4)-1(4,1)+4(4,2)+5(4,2,1)-1(4,3)+2(4,3,1)\\
s_{17} & = & (1)+6(2)+5(2,1)+15(3,1)+11(3,2)+3(3,2,1) 
           -12(4)-2(4,1)+4(4,2)+5(4,2,1)+4(4,3,1)\\
s_{18} & = & (1)+7(2)+5(2,1)+(3)+17(3,1)+11(3,2)+4(3,2,1) 
           -12(4)-2(4,1)+2(4,2)+7(4,2,1)-2(4,3)+4(4,3,1)+(4,3,2)\\
s_{19} & = & (1)+6(2)+6(2,1)+19(3,1)+14(3,2)+4(3,2,1) 
           -12(4)-3(4,1)+4(4,2)+10(4,2,1)-3(4,3)+5(4,3,1)+(4,3,2)\\
s_{20} & = & (1)+7(2)+6(2,1)+21(3,1)+15(3,2)+5(3,2,1) 
           -17(4)-3(4,1)+5(4,2)+11(4,2,1)-3(4,3)+8(4,3,1)+(4,3,2)\\
s_{21} & = & (1)+8(2)+6(2,1)+(3)+23(3,1)+16(3,2)+6(3,2,1) 
           -20(4)-3(4,1)+4(4,2)+13(4,2,1)-4(4,3)+10(4,3,1)+2(4,3,2)\\
s_{22} & = & (1)+7(2)+7(2,1)+26(3,1)+19(3,2)+6(3,2,1) 
           -20(4)-4(4,1)+7(4,2)+17(4,2,1)-6(4,3)+11(4,3,1)+3(4,3,2)\\
s_{23} & = & (1)+8(2)+7(2,1)+28(3,1)+21(3,2)+7(3,2,1)-28(4) 
           -4(4,1)+10(4,2)+18(4,2,1)-5(4,3)+16(4,3,1)+3(4,3,2)\\
s_{24} & = & (1)+9(2)+7(2,1)+(3)+30(3,1)+23(3,2)+8(3,2,1)-31(4)
           -6(4,1)+11(4,2)+20(4,2,1)-5(4,3)+19(4,3,1)+5(4,3,2)\\
s_{25} & = & (1)+8(2)+8(2,1)+32(3,1)+25(3,2)+9(3,2,1)-28(4) 
           -9(4,1)+10(4,2)+24(4,2,1)-7(4,3)+20(4,3,1)+8(4,3,2)\\
s_{26} & = & (1)+9(2)+8(2,1)+35(3,1)+27(3,2)+10(3,2,1)-36(4) 
           -8(4,1)+10(4,2)+26(4,2,1)-7(4,3)+24(4,3,1)+7(4,3,2)+(4,3,2,1)
\end{eqnarray*}
}

With mod $2$ reduction (the vectors are denoted by $sq_i)$:
{\scriptsize
\begin{eqnarray*}
sq_{0} & = & 1(0)\\
sq_{1} & = & (1)\\
sq_{2} & = & (1)+(2)\\
sq_{3} & = & (1)+(3)\\
sq_{4} & = & (1)+(2)+(2,1)+(4)\\
sq_{5} & = & (1)+(2,1)+(3,1)\\
sq_{6} & = & (1)+(2)+(2,1)+(3)+(4,1)\\
sq_{7} & = & (1)+(3,2)\\
sq_{8} & = & (1)+(2)+(3,2)+(4)+(4,1)+(4,2)\\
sq_{9} & = & (1)+(3)+(3,1)+(4,1)\\
sq_{10} & = & (1)+(2)+(2,1)+(3,1)+(4,1)+(4,2)+(4,3)\\
sq_{11} & = & (1)+(2,1)+(3,2)+(3,2,1)+(4,3)\\
sq_{12} & = & (1)+(2)+(2,1)+(3)+(3,2,1)+(4)+(4,2)+(4,2,1)\\
sq_{13} & = & (1)+(3,1)+(3,2,1)+(4,1)\\
sq_{14} & = & (1)+(2)+(4,1)+(4,3,1)\\
sq_{15} & = & (1)+(3)+(3,1)+(3,2)+(4,3,1)\\
sq_{16} & = & (1)+(2)+(2,1)+(4)+(4,1)+(4,2,1)+(4,3)\\
sq_{17} & = & (1)+(2,1)+(3,1)+(3,2)+(3,2,1)+(4,2,1)\\
sq_{18} & = & (1)+(2)+(2,1)+(3)+(3,1)+(3,2)+(4,2,1)+(4,3,2)\\
sq_{19} & = & (1)+(3,1)+(4,1)+(4,3)+(4,3,1)+(4,3,2)\\
sq_{20} & = & (1)+(2)+(3,1)+(3,2)+(3,2,1)+(4)+(4,1) 
 +(4,2)+(4,2,1)+(4,3)+(4,3,2)\\
sq_{21} & = & (1)+(3)+(3,1)+(4,1)+(4,2,1)\\
sq_{22} & = & (1)+(2)+(2,1)+(3,2)+(4,2)+(4,2,1)+(4,3,1)+(4,3,2)\\
sq_{23} & = & (1)+(2,1)+(3,2)+(3,2,1)+(4,3)+(4,3,2)\\
sq_{24} & = & (1)+(2)+(2,1)+(3)+(3,2)+(4)+(4,2)+(4,3)+(4,3,1)+(4,3,2)\\
sq_{25} & = & (1)+(3,2)+(3,2,1)+(4,1)+(4,3)\\
sq_{26} & = & (1)+(2)+(3,1)+(3,2)+(4,3)+(4,3,2)+(4,3,2,1)\\
\end{eqnarray*}
}

\subsection{Transposed matrices of the endomorphism $[\varrho] \mapsto  \sum_{0 \leq i \leq n}\, \, [\varrho \otimes \Lambda^i]$ of $R_{\F2}(\GL_n(\F2))$
}

\

 For $n=0$, the basis being   $\{(0)\}$,

$$\tau_0=\left( \begin {array}{c} 1\end {array} \right) 
$$

For $n=1$, the basis being   $\{(1)\}$,

$$\tau_1=\left( \begin {array}{c} 2\end {array} \right) 
$$

For $n=2$, the basis being   $\{(2), (2,1)\}$,

$$ \tau_2=\left( \begin {array}{cc} 2&2\\ \noalign{\medskip}1&3\end {array}
 \right) 
$$

For $n=3$, the basis being   $\{(3), (3,1), (3,2), (3,2,1)\}$,

$$\tau_3=\left( \begin {array}{rrrr} 2&0&0&0\\ \noalign{\medskip}1&3&2&6
\\ \noalign{\medskip}1&2&3&6\\ \noalign{\medskip}0&1&1&4\end {array}
 \right)
$$

For $n=4$, the basis being   $\{(4), (4,1), (4,2), (4,2,1),(4,3), (4,3,1), (4,3,2),(4,3,2,1)\}$,

$$\tau_4=\left( \begin {array}{rrrrrrrr} 2&0&-2&-4&0&0&-4&-16
\\ \noalign{\medskip}1&3&0&-2&0&0&0&-4\\ \noalign{\medskip}1&2&3&2&2&0
&2&8\\ \noalign{\medskip}0&1&1&4&0&2&4&14\\ \noalign{\medskip}1&0&0&0&
3&0&-2&-4\\ \noalign{\medskip}0&1&2&6&1&4&6&24\\ \noalign{\medskip}0&0
&1&4&1&2&4&14\\ \noalign{\medskip}0&0&0&1&0&1&1&5\end {array} \right)
$$

   For $n=5$, the basis being 
   $\{(5)$, $(5,1)$, $(5,2)$, $(5,2,1)$, $(5,3)$, $(5,3,1)$, $(5,3,2)$, $(5,3,2,1)$, 
      $(5,4),$ $(5,4,1)$, $(5,4,2)$, $(5,4,2,1)$,$(5,4,3)$, $(5,4,3,1)$, $(5,4,3,2)$,$(5,4,3,2,1)\}$,
   $$
\tau_5 =
 \left( \begin{smallmatrix*}[r] 2&0&0&0&0&0&0&0&0&0&0&0&0&0&0
&0\\ \noalign{\medskip}1&3&0&0&-2&-2&-4&-12&0&0&0&0&-4&-4&-8&-56
\\ \noalign{\medskip}1&2&3&2&0&-4&-8&-28&0&0&0&0&0&-8&-32&-152
\\ \noalign{\medskip}0&1&1&4&0&0&-2&-10&0&0&0&0&0&0&-4&-36
\\ \noalign{\medskip}1&0&0&0&3&0&-8&-32&2&0&-4&-8&2&0&-28&-152
\\ \noalign{\medskip}0&1&2&6&1&4&2&-10&0&2&0&-4&0&2&8&-4
\\ \noalign{\medskip}0&0&1&4&1&2&4&6&0&0&2&0&4&0&6&32
\\ \noalign{\medskip}0&0&0&1&0&1&1&5&0&0&0&2&0&4&8&30
\\ \noalign{\medskip}1&0&-2&-4&0&0&-4&-8&3&0&-2&-4&0&0&-12&-56
\\ \noalign{\medskip}0&1&0&-2&0&0&0&-4&1&4&0&-2&-2&-2&-4&-24
\\ \noalign{\medskip}0&0&1&0&2&0&2&8&1&2&4&2&6&-4&-10&-4
\\ \noalign{\medskip}0&0&0&1&0&2&4&14&0&1&1&5&0&6&18&80
\\ \noalign{\medskip}0&0&0&0&1&0&-2&-4&1&0&0&0&4&0&-10&-36
\\ \noalign{\medskip}0&0&0&0&0&1&4&18&0&1&2&6&1&5&14&80
\\ \noalign{\medskip}0&0&0&0&0&0&1&8&0&0&1&4&1&2&5&30
\\ \noalign{\medskip}0&0&0&0&0&0&0&1&0&0&0&1&0&1&1&6
\end{smallmatrix*}
 \right) 
$$

For $n=6$, the basis being generated in the same way: 
$$ \tau_6 =
\left(\begin{smallmatrix*}[r] 2&0&0&0&2&0&0
&0&0&0&0&0&0&0&0&0&0&0&0&0&0&0&0&0&0&0&0&0&0&0&0&0
\\ \noalign{\medskip}1&3&0&0&0&2&4&12&0&0&0&0&0&0&-8&-72&0&0&0&0&0&0&0
&0&0&0&0&0&0&0&-16&-272\\ \noalign{\medskip}1&2&3&2&0&0&10&40&-2&0&2&8
&-4&0&4&-128&0&0&0&0&0&0&0&-32&-4&0&-4&0&-16&16&80&-80
\\ \noalign{\medskip}0&1&1&4&0&0&0&10&0&-2&0&2&0&-4&-8&-20&0&0&0&0&0&0
&0&0&0&-4&0&-4&8&-8&-40&-184\\ \noalign{\medskip}1&0&0&0&3&0&0&0&0&0&0
&0&0&0&0&0&0&0&0&0&0&0&0&0&0&0&0&0&0&0&0&0\\ \noalign{\medskip}0&1&2&6
&1&4&2&6&0&0&-4&-4&-8&-8&-28&-124&0&0&0&0&0&0&0&0&0&0&-8&-8&-32&0&-24&
-568\\ \noalign{\medskip}0&0&1&4&1&2&4&6&0&0&0&-8&-2&-16&-30&-128&0&0&0
&0&0&0&0&0&0&0&0&-16&-4&-64&-196&-1008\\ \noalign{\medskip}0&0&0&1&0&1
&1&5&0&0&0&0&0&-2&-8&-38&0&0&0&0&0&0&0&0&0&0&0&0&0&-4&-32&-228
\\ \noalign{\medskip}1&0&-2&-4&0&0&-4&-16&3&0&2&-4&10&0&4&80&2&0&0&0&0
&0&0&16&2&0&8&0&40&-32&-128&-80\\ \noalign{\medskip}0&1&0&-2&0&0&0&-4&
1&4&0&2&-8&2&4&-48&0&2&0&0&-4&-4&-8&-24&0&2&0&8&-28&12&136&280
\\ \noalign{\medskip}0&0&1&0&2&0&2&8&1&2&4&2&2&-16&-52&-272&0&0&2&0&0&
-8&-16&-56&0&0&2&0&8&-56&-272&-1608\\ \noalign{\medskip}0&0&0&1&0&2&4&
14&0&1&1&5&0&2&-12&-120&0&0&0&2&0&0&-4&-20&0&0&0&2&0&8&-12&-412
\\ \noalign{\medskip}0&0&0&0&1&0&-2&-4&1&0&0&0&4&0&-30&-196&0&0&0&0&2&0
&-16&-64&4&0&-8&-16&6&0&-128&-1008\\ \noalign{\medskip}0&0&0&0&0&1&4&
18&0&1&2&6&1&5&6&-68&0&0&0&0&0&2&0&-32&0&4&0&-8&0&6&32&-80
\\ \noalign{\medskip}0&0&0&0&0&0&1&8&0&0&1&4&1&2&5&14&0&0&0&0&0&0&2&0&0
&0&4&0&8&0&14&96\\ \noalign{\medskip}0&0&0&0&0&0&0&1&0&0&0&1&0&1&1&6&0
&0&0&0&0&0&0&2&0&0&0&4&0&8&16&62\\ \noalign{\medskip}1&0&0&0&0&0&0&0&0
&0&0&0&4&0&-8&-16&3&0&0&0&2&0&0&0&0&0&0&0&12&0&-72&-272
\\ \noalign{\medskip}0&1&0&0&-2&-2&-4&-12&0&0&0&0&-4&0&16&88&1&4&0&0&-
2&0&0&0&0&0&0&0&-12&0&88&544\\ \noalign{\medskip}0&0&1&0&0&-4&-8&-28&0
&0&0&0&0&-8&4&136&1&2&4&2&0&-4&2&12&-2&0&2&8&-4&-24&-48&280
\\ \noalign{\medskip}0&0&0&1&0&0&-2&-10&0&0&0&0&0&0&-4&0&0&1&1&5&0&0&-
2&0&0&-2&0&2&0&-4&-20&-84\\ \noalign{\medskip}0&0&0&0&1&0&-8&-32&2&0&-
4&-8&2&0&-28&-24&1&0&0&0&4&0&-8&0&6&0&-4&-8&6&0&-124&-568
\\ \noalign{\medskip}0&0&0&0&0&1&0&-16&0&2&0&-4&0&2&8&-4&0&1&2&6&1&5&2
&-10&0&6&-4&-8&-16&-10&-4&-48\\ \noalign{\medskip}0&0&0&0&0&0&1&0&0&0&
2&0&4&0&6&32&0&0&1&4&1&2&5&6&0&0&6&-8&18&-32&-68&-80
\\ \noalign{\medskip}0&0&0&0&0&0&0&1&0&0&0&2&0&4&8&30&0&0&0&1&0&1&1&6&0
&0&0&6&0&18&56&260\\ \noalign{\medskip}0&0&0&0&0&0&0&8&1&0&0&0&0&0&-8&
-40&1&0&-2&-4&0&0&-4&-8&4&0&2&-4&10&0&-20&-184\\ \noalign{\medskip}0&0
&0&0&0&0&0&0&0&1&0&0&-2&-2&-4&-20&0&1&0&-2&0&0&0&-4&1&5&0&2&-10&0&0&-
84\\ \noalign{\medskip}0&0&0&0&0&0&0&0&0&0&1&0&4&-4&-12&-12&0&0&1&0&2&0
&2&8&1&2&5&2&14&-20&-120&-412\\ \noalign{\medskip}0&0&0&0&0&0&0&0&0&0&0
&1&0&4&14&66&0&0&0&1&0&2&4&14&0&1&1&6&0&14&66&402\\ \noalign{\medskip}0
&0&0&0&0&0&0&0&0&0&0&0&1&0&-8&-32&0&0&0&0&1&0&-2&-4&1&0&0&0&5&0&-38&-
228\\ \noalign{\medskip}0&0&0&0&0&0&0&0&0&0&0&0&0&1&8&56&0&0&0&0&0&1&4
&18&0&1&2&6&1&6&30&260\\ \noalign{\medskip}0&0&0&0&0&0&0&0&0&0&0&0&0&0
&1&16&0&0&0&0&0&0&1&8&0&0&1&4&1&2&6&62\\ \noalign{\medskip}0&0&0&0&0&0
&0&0&0&0&0&0&0&0&0&1&0&0&0&0&0&0&0&1&0&0&0&1&0&1&1&7 \end{smallmatrix*}\right)$$

\end{document}